\documentstyle[amstex,amssymb,12pt,leqno,xcolor]{article}
\newcommand {\R} {{\mathbb R}}
\newcommand {\C} {{\mathbb C}}
\newcommand {\N} {{\mathbb N}}

\renewcommand {\H}{{\mathbb H}}

\newcommand{\essinf}{{{\rm ess\,inf\,}}}

\newtheorem{proposition}{Proposition}[section]
\newtheorem{theorem}{Theorem}

\newtheorem{lemma}[proposition]{Lemma}

\def    \R      {{\mathbb R}}

\newcommand{\grad}{{\rm grad}}

\def\kasten{\hfill\null\nobreak\hfill \hbox{\vrule\vbox{\hrule width
6pt\vskip6pt\hrule}\vrule} \par\smallskip}

\begin{document}
\begin{center}
{\Large {\bf Reflection of Willmore surfaces\\
\vspace*{3mm}with free boundaries}}
\vspace*{5mm}\\
{\large {\sc Ernst Kuwert}} \& {\large {\sc Tobias Lamm}}\vspace*{5mm}\\
\end{center}
\begin{abstract}
We study immersed surfaces in $\R^3$ which are critical points of the Willmore functional under boundary constraints. The two cases considered are when the surface meets a plane orthogonally along the boundary, and when 
the boundary is contained in a line. In both cases 
we derive weak forms  of the resulting free boundary conditions and prove regularity by reflection.
\end{abstract}

\section{Introduction}

This note is concerned with Willmore surfaces under free boundary 
conditions. Let $\Sigma$ be an oriented, two-dimensional manifold 
with boundary $\partial \Sigma$. For a smooth immersion $f:\Sigma \to \R^3$
the Willmore functional is defined by 
\begin{equation}
\label{eqwillmorefunctional}
{\cal W}(f) = \frac{1}{4} \int_\Sigma H^2\,d\mu_g.
\end{equation}
Here $\mu_g$ is the measure associated to the induced Riemannian
metric $g$, and $H$ is the mean curvature with respect to the unit normal 
$\nu:\Sigma \to {\mathbb S}^2$. We denote by $h$ the second fundamental
form of $f$ and put $h^\circ = h - \frac{1}{2} H g$.
Let $f(\cdot,t)$ be a smooth variation with velocity field
$\phi = \varphi \nu + Df \cdot \xi$. The first variation 
is, see \cite{Nit93} and \cite{AK14},
\begin{equation}
\label{eqwillmorefirst}
\frac{d}{dt} {\cal W}(f(\cdot,t))|_{t=0} = \frac{1}{2} \int_\Sigma W(f) \varphi\,d\mu_g
+ \frac{1}{2} \int_{\partial \Sigma} \omega(\eta)\,ds_g ,
\end{equation}
where $\eta$ denotes the interior unit normal along $\partial \Sigma$ with 
respect to $g$, and 
\begin{eqnarray} 
\label{eqwillmoreoperator} 
W(f) & = & \Delta_g H + |h^\circ|^2 H,\\ 
\label{eqwillmoreboundaryform}
\omega(\eta) & = & 
\varphi \frac{\partial H}{\partial \eta} 
- \frac{\partial \varphi}{\partial \eta} H
- \frac{1}{2} H^2 g(\xi,\eta).
\end{eqnarray}
In this paper we address two free boundary situations. First,
for a given support surface $S$ we consider the class of smooth
immersions $\widetilde{{\cal M}}(S)$ which meet $S$ orthogonally along the boundary $\partial \Sigma$.
The admissible variations $\phi \in T_f{\widetilde{\cal M}}(S)$ are 
characterized by the equations
\begin{equation}
\label{eqadmiissible}
\frac{\partial \varphi}{\partial \eta} + h^S(\nu,\nu) \varphi = 0,\quad
g(\xi,\eta) = 0.
\end{equation}
For immersions $f$ which are critical in that class, Alessandroni and the first author computed
in \cite{AK14} the following natural boundary condition
\begin{equation}
\label{eqfreeboundarycondition}
\frac{\partial H}{\partial \eta} + h^S(\nu,\nu) H = 0  \quad 
\mbox{ along }\partial \Sigma,
\end{equation}
where $h^S$ is the second fundamental form of $S \subset \R^3$.

Secondly, we consider immersions which are confined to a 
given support curve $\Gamma \subset \R^3$ along the boundary, but 
without prescribing the tangent plane along $\partial \Sigma$. In 
this case the critical immersions in the corresponding class 
satisfy 
\begin{equation}
\label{eqnavierboundary}
H = 0 \quad \mbox{ along } \partial \Sigma.
\end{equation}
This is referred to as Navier boundary condition in \cite{DG09}.
Our note deals with the cases of a plane $S$ and a line $\Gamma$, 
proving reflection principles in both situations. For the plane the result 
is already due to J.C.C. Nitsche, assuming $C^{4,\nu}$ regularity up to
the boundary \cite{Nit91}. We extend Nitsche's theorem to a weak setting,
and further adapt our arguments to the case of a line $\Gamma$.\\
\\
The class of $W^{2,2}$ (Lipschitz) immersions on an open set 
$U \subset \R^2$ is defined by 
\begin{equation}
\label{eqweakimmersion}
W^{2,2}_{{\rm imm}}(U,\R^3) = \{f \in W^{2,2} \cap W^{1,\infty}(U,\R^3): 
\essinf (\det g) > 0\}.  
\end{equation}
This is an open subset of $W^{2,2}\cap W^{1,\infty}(U,\R^3)$. The
$W^{2,2}$ conformal immersions $W^{2,2}_{{\rm conf}}(U,\R^3)$ introduced
in \cite{KL12} are exactly those $f \in  W^{2,2}_{{\rm imm}}(U,\R^3)$
satisfying $g_{11} = g_{22}$ and $g_{12} = 0$.
In other words, we have $g_{ij} = e^{2u} \delta_{ij}$
where $u \in W^{1,2} \cap L^\infty(U)$. In the following we let
$Q = (-\pi,\pi) \times (-1,1)$, $I = (-\pi,\pi) \times \{0\} =: (-\pi,\pi)$
and $Q_{\pm} = \{(x,y) \in Q: {\pm} y > 0\}$. The precise choice
of $Q$ will be convenient in the appendix.

{\bf Theorem } {\em Let ${\cal M}$ be the class of 
$f \in W^{2,2}_{{\rm imm}}(Q_{+},\R^3)$ satisfying 
the constraints
\begin{eqnarray}
\label{eqconstraint1}
f_3 & = & 0 \quad \mbox{ along }I,\\ 
\label{eqconstraint2}
\langle \nu , e_3 \rangle & = & 0 \quad 
\mbox{ along }I \mbox{ almost everywhere.}
\end{eqnarray}
Assume that $f \in W^{2,2}_{{\rm conf}}(Q_{+},\R^3)$ is a critical 
point of the Willmore energy for variations in ${\cal M}$ with 
compact support in $Q_{+} \cup I$. Then extending $f$ to $Q$ by
$$
f:Q_{-} \to \R^3,\,f(x,y) = (f_1(x,-y),f_2(x,-y),-f_3(x,-y)),
$$
yields a smooth Willmore immersion $f \in C^\infty(Q,\R^3)$.}\\
\\
The key tool of the proof is the interior regularity
theorem of Rivi\`{e}re \cite{Riv08}, which directly 
implies regularity on $Q_{+}$ (and on $Q_{-}$).
We show that $f$ is of class $W^{2,2}_{{\rm conf}}(Q,\R^3)$ 
and solves the weak Willmore equation on all of $Q$, so 
that \cite{Riv08} applies. 
For the second free boundary problem where $f$ is critical
among immersions mapping $I$ into a line the arguments 
are quite analogous. The extension by reflection across 
the line again gives a smooth Willmore immersion, see
Theorem \ref{thm2} in Section \ref{section3}.\\
\\
The notion of critical point used in Theorem 1 is that 
the first variation $\delta {\cal W}(f,\phi)$ vanishes
for all {\em admissible} vector fields $\phi$. Here
admissible means that $\phi$ is formally a tangent
vector to ${\cal M}$ at $f$, i.e. it satisfies
the equations obtained by linearizing the constraints
(\ref{eqconstraint1}) and (\ref{eqconstraint2}) 
at the immersion $f \in {\cal M}$ (clearly only 
(\ref{eqconstraint2}) is non-linear). 
We prove in Section \ref{section2} that locally any 
admissible vector field $\phi$ is indeed the 
tangent vector of a curve in ${\cal M}$ at $f$.
Therefore our regularity result would apply, for 
example, to show the regularity of minimizers.\\
\\
In the case where the curve and the tangent plane are
prescribed along the boundary, there are substantial 
existence and regularity results for minimizers by
Sch\"atzle \cite{Sch10} and Da Lio, Palmurella and 
Rivi\`{e}re \cite{DPR18}. It is clearly of interest 
to develop an analogous theory for the free boundary 
problems in the case of curved supporting surfaces 
or curves.\\ 
\\
Due to its conformal invariance it is also interesting to study the corresponding problem for the functional 
$$
{\cal T}(f) = \frac{1}{2} \int_\Sigma |h^\circ|^2\,d\mu_g.
$$
In section $5$ we calculate the first variation of ${\cal T}$ to be
\begin{equation}
\label{eqsfffirst}
\frac{d}{dt} {\cal T}(f(\cdot,t))|_{t=0} = \frac{1}{2} \int_\Sigma W(f) \varphi\,d\mu_g
+  \int_{\partial \Sigma} \alpha(\eta)\,ds_g ,
\end{equation}
where $W(f)$ and the smooth variation $f(\cdot,t)$ are as above and 
\begin{eqnarray}  
\label{eqsffboundaryform}
\alpha(\eta) & = & 
\frac12 \varphi \frac{\partial H}{\partial \eta} 
- h^\circ(\grad \varphi,\eta)-\frac12 |h^\circ|^2g(\xi,\eta).
\end{eqnarray}
We consider again the free boundary problem when the immersions meet $S$ orthogonally along the boundary $\partial \Sigma$. 
In the case of the support surface $S$ being a plane the boundary condition is again given by
$$
\frac{\partial H}{\partial \eta}=0 \ \ \ \text{along}\ \ \ \partial \Sigma.
$$
It follows immediately from the above that our main theorem directly extends to critical points of ${\cal T}$ in this situation.

Assume that $f:D\to \R^3$ is a conformally immersed disk which is critical in $\cal{M}$. We have shown that reflection 
at $\R^2$ extends the surface to a Willmore immersion 
$f:\widehat{\C} \to \R^3$. Bryant proved that there exists 
a round sphere, such that under the associated inversion
one obtains a complete minimal immersion of finite 
total curvature, with a finite number of flat ends \cite{Bry84}. 
In fact Bryant's theory includes the case of branched 
immersions under various assumptions, see \cite{LN}, \cite{Riv}. 
We claim that the center of a sphere with these properties
is unique. Otherwise, by scaling and rotating, we can assume 
that $f_{\pm}: = I_{\pm} \circ f$ are minimal surfaces of 
that type, where $I_{\pm}$ are the inversions at the spheres 
of radiius $\sqrt{2}$ around $\pm e_3$. It follows that 
$f_{+} = I \circ f_{-}$ where $I:= I_{+} \circ I_{-}$ 
is the inversion $I(x) = \frac{x}{|x|^2}$. Now in general 
we have the relation, see (2.27) in \cite{BK03},
$$
\frac{1}{4} |\vec{H}_{+}|^2\,d\mu_{+} =
\frac{1}{4} |\vec{H}_{-}|^2\,d\mu_{-} + 
\big(\Delta_{g_{-}} \log |f_{-}|^2\big)\,d\mu_{-}.
$$
But in our case $H_{+} = H_{-} = 0$, and 
hence we get, away from finitely many points,
$$
0 = \Delta_{g_{-}}\,\log |f_{-}|^2  = 
\frac{4 |f_{-}^\perp|^2}{|f_{-}|^4}.
$$
It follows that $f_{-}^\perp$ vanishes, which means that 
$f_{\pm}$ are conical about the origin. 
But there is no smooth minimal cone in $\R^3$ 
except the plane. Thus unless our initial surface is 
a round half-sphere, the center of inversion lies on $\R^2$ 
and $f$ is an inverted minimal surface with $\R^2$ symmetry.

In Bryant's list, the first example is the Morin surface, 
which is symmetric even under a rotation by $\frac{\pi}{2}$. It is the inversion of (see e.g. \cite{kusner}) 
$$
f:{\mathbb C}^2\backslash \{p_1,\ldots,p_4\} \to \R^3, \, f(w)=\Re \Big( \frac{i(w^3-w),(w^3+w),\frac{i}{2}(w^4+1)}{w^4+2\sqrt{3}w^2-1} \Big),
$$
where the $p_i$, $1\le i\le 4$, are the zeros of the function in the demoninator of this expression. It is easy to see that the image of $f$ contains the $y$-axis and hence the Morin surface is invariant under reflections at this line and thus we found an example for the case in which we confine the immersion to a given support line. Additionally, it follows that the conjugate surface $f^\star$ is invariant under reflections at the $xz$-plane. This inversion of $f^\star$ thus yields an example for the other situation considered in this paper.

Another example for the reflection at a plane is obtained from the catenoid with 
conformal parametrization
$$
f(s,\theta) = 2\,(\cosh s \cos \theta,\cosh s \sin \theta, s)
\quad \mbox{ for }s,\theta \in \R.
$$ 
Inverting the catenoid yields a bounded surface which has a 
double point at the origin with horizontal tangent plane,
corresponding to $s \to \pm \infty$. 
Substituting $s+i \theta = \log w$ where $w = \varrho e^{i\theta}$, 
we obtain for the induced metric
$$
g_{ij} = \frac{1+2\varrho^2 + \varrho^4}
{(1+ 2\varrho^2 + \varrho^4 + 4\varrho^2 \log^2 \varrho)^2}\,\delta_{ij} 
= \big(1-4\varrho^2 \log^2 \varrho + {\mathcal O}(\varrho^2)\big)\,\delta_{ij}
\quad \mbox{ as }\varrho \to 0.
$$
This parametrization is $W^{2,2}$-conformally 
immersed near $w = 0$. Now restricting to 
$-\frac{\pi}{2} \leq \theta \leq \frac{\pi}{2}$ gives 
a surface which meets the vertical plane $x_1 = 0$ 
orthogonally. Moreover for $\theta = \pm \frac{\pi}{2}$
the free boundary condition (\ref{eqfreeboundarycondition})
holds, namely we have $\partial_\eta H = 0$ by rotational
symmetry. But for this we have to exclude the origin, 
in fact the mean curvature vector has a singular 
expansion \cite{KS04}
$$
\vec{H} = - (\log \varrho)\,e_3 + {\cal O}(1) \quad
\mbox{ as } \varrho \to 0.
$$
If the first variation would vanish for all admissible 
vector fields, then we could apply our main Theorem to conclude that reflection at the plane $x_1=0$ produces a smooth extension.
Thus the example only matches our assumptions away 
from the two singular points. Therefore, the regularity result does not follow if the boundary condition is only satisfied away from a point.

\section{Reflection at a plane}
\label{section1}

Let us start by recalling that for $u \in W^{1,2}(Q_{\pm})$ we 
have the upper and lower traces $u_{\pm} \in L^2(I)$. We note 
that $u_{-} = \pm u_{+}$ for $u$ even resp. odd. 

\begin{lemma} \label{lemmaw12reflection} The following holds 
for any $u \in W^{1,2}(Q_{\pm})$:
\begin{itemize}
\item[{\rm (1)}] $\partial_x u$ is the weak derivative on $Q$.
\item[{\rm (2)}] $\partial_y u$ is the weak derivative on $Q$
                 if and only if $u_{\pm}$ coincide.
\end{itemize}
\end{lemma}

{\em Proof. } We have $u(\cdot,y) \in W^{1,2}(I)$ for almost every $y \in I$.
For statement (1) we compute for $\varphi \in C^\infty_c(Q)$ using Fubini
\begin{eqnarray*}
\int_Q u \partial_x \varphi\,dxdy & = & 
\int_{-1}^1 \Big(\int_{-\pi}^{\pi} u(x,y) \partial_x \varphi(x,y)\,dx\Big)\,dy\\
& = & - \int_{-1}^1 \Big(\int_{-\pi}^{\pi} \partial_x u(x,y) \varphi(x,y)\,dx\Big)\,dy\\
& = & \int_Q (\partial_x u) \varphi\,dxdy.
\end{eqnarray*}
For statement (2) we use
\begin{eqnarray*}
\int_Q u \partial_y \varphi & = &
\int_{Q_{+}} u \partial_y \varphi + \int_{Q_{-}} u \partial_y \varphi\\
& = & -\int_I u_{+} \varphi\,dx - \int_{Q_{+}} (\partial_y u) \varphi
+ \int_I u_{-} \varphi \,dx - \int_{Q_{-}} (\partial_y u) \varphi\\ 
& = & \int_I (u_{-}-u_{+}) \varphi\,dx - \int_Q (\partial_y u) \varphi.
\end{eqnarray*}
\kasten

By the Sobolev embedding theorem, a function $u \in W^{2,2}(Q_{\pm})$ has a 
representative in $C^{0,\alpha}(\overline{Q}_{\pm})$, for any $\alpha \in [0,1)$. 
In particular the traces are given by continuous extension. Furthermore the 
derivatives have traces $(\partial_i u)_{\pm} \in L^2(I)$. In the following 
we need the concept of $W^{2,2}$ conformal immersions, see \cite{KL12}. 

\begin{lemma} \label{lemmaw22conformal} Let $f:Q_{+} \to \R^3$ be a $W^{2,2}$ 
conformal immersion meeting the horizontal plane $\R^2$ orthogonally along
$I$, that is
\begin{eqnarray}
\label{eqboundary1}
f_3 & = & 0 \quad \mbox{ on }I, \\
\label{eqboundary2} 
\langle \nu_{+}, e_3 \rangle & = & 0 \quad \mbox{ a.e. on $I$}\quad \mbox{ where } 
\nu = \frac{\partial_x f \times \partial_y f}{|\partial_x f \times \partial_y f|}.
\end{eqnarray}
Let $f$ be extended to $Q$ by reflection at $\R^2$, that is 
\begin{equation}
\label{eqextend}
f_i(x,y)  = \,\, f_i(x,-y) \,\mbox{ for } i = 1,2, \quad
f_3(x,y) = - f_3(x,-y). 
\end{equation}
Then $f:Q \to \R^3$ is a $W^{2,2}$ conformal immersion. 
\end{lemma}

{\em Proof. }We first note that 
$\partial_x f \times \partial_y f \in W^{1,2} \cap L^\infty(Q_{+},\R^3)$ and
$$
\essinf |\partial_x f \times \partial_y f| \geq e^{2\lambda} > 0 \quad 
\mbox{ for some }\lambda \in \R.
$$
The Sobolev chain rule yields that $\nu \in W^{1,2}(Q_{+},\R^3)$ 
so that the trace $\nu_{+} \in L^\infty(I)$ is defined.
The functions $f_i$, $i = 1,2$, are even, while $f_3$ is odd 
and vanishes on $I$ by assumption (\ref{eqboundary1}).
Thus Lemma \ref{lemmaw12reflection} yields directly $f \in W^{1,2}(Q,\R^3)$.
To see that  $f$ is actually of class $W^{2,2}(Q,\R^3)$ we show that
\begin{eqnarray}
\label{eqboundary3}
(\partial_x f_3)_{+} & = & 0  \quad \mbox{ on }I,\\ 
\label{eqboundary4}
(\partial_y f_i)_{+} & = & 0  \quad \mbox{ on }I \quad \mbox{ for } i = 1,2.
\end{eqnarray}
These are the derivatives which are odd, the others are even. 
For given $\varphi \in C^\infty_c(Q_{+} \cup I)$ we calculate 
using the divergence theorem and partial integration w.r.t. $\partial_x$
\begin{eqnarray*}
\int_I (\partial_x f_3)_{+} \varphi\,dx & = & 
- \int_{Q_{+}} \partial_y (\partial_x f_3)\, \varphi\,dxdy
- \int_{Q_{+}} (\partial_ x f_3)\, \partial_y \varphi\,dxdy\\
& = & - \int_{Q_{+}} \partial_x (\partial_y f_3)\, \varphi\,dxdy
+ \int_{Q_{+}} f_3\, \partial_x (\partial_y \varphi)\,dxdy\\
& = & \int_{Q_{+}} (\partial_y f_3)\, \partial_x \varphi\,dxdy
+ \int_{Q_{+}} f_3\, \partial_y (\partial_x \varphi)\,dxdy\\
& = & - \int_I f_3\, \partial_x \varphi\,dx\\ 
& = & 0.
\end{eqnarray*}
Since $\varphi$ is arbitrary we get (\ref{eqboundary3}), more 
precisely we have shown that
$$
(\partial_x f_3)_{+} = \partial_x (f_3|_I) = 0.
$$
Next we claim that almost everywhere on $I$ we have the equations
\begin{eqnarray}
\label{eqconformboundary}
\big|(\partial_x f)_{+}|^2 - \big|(\partial_y f)_{+}|^2 & = &
\big\langle (\partial_x f)_{+},(\partial_y f)_{+} \big\rangle = 0,\\
\label{eqimmersedboundary}
\big|(\partial_x f)_{+} \times (\partial_y f)_{+}\big| & \geq & e^{2\lambda}.
\end{eqnarray}
To see this, we use that for any $u \in W^{1,2}(Q_{+})$ there exists a 
sequence $y_k \searrow 0$ such that
$$
u(\cdot,y_k) \to u_{+} \quad \mbox{ almost everywhere and in }L^2(I).
$$
More precisely, there is a null set $N \subset (0,1)$ such that 
any sequence $y_k \searrow 0$ with $y_k \notin N$ has a subsequence 
with this property. Therefore we can chose $y_k \searrow 0$ with 
$$
\partial_x f(\cdot,y_k) \to (\partial_x f)_{+},\quad 
\partial_y f(\cdot,y_k) \to (\partial_y f)_{+},\quad
\nu(\cdot,y_k) \to \nu_{+}.
$$
We can further assume that almost everywhere on $I$
$$
|\partial_x f(\cdot,y_k)|^2 - |\partial_y f(\cdot,y_k)|^2 = 
\langle \partial_x f(\cdot,y_k), \partial_y f(\cdot,y_k) \rangle  =0,
\quad \mbox{ and }
$$
$$
|\partial_x f(\cdot,y_k) \times \partial_y f(\cdot,y_k)| \geq e^{2\lambda}.
$$
Passing to limits proves (\ref{eqconformboundary}) and (\ref{eqimmersedboundary}). 
We now verify (\ref{eqboundary4}) for a.e. $(x,0) \in I$.
We may assume $(\partial_y f)_{+}(x,0) \neq 0$, hence 
(\ref{eqconformboundary}) and (\ref{eqboundary3}) imply 
$(\partial_x f)_{+}(x,0) \in \R^2 \backslash \{0\}$. Now 
$$
\nu_{+} = 
\frac{(\partial_x f)_{+} \times (\partial_y f)_{+}}{|(\partial_x f)_{+} 
\times (\partial_y f)_{+}|} 
\quad \mbox{ almost everywhere on }I.
$$
By (\ref{eqboundary2}) we get that $(\partial_y f)_{+}(x,0)$ is a linear
combination of $(\partial_x f)_{+}(x,0)$ and $e_3$. But then by 
(\ref{eqconformboundary}) the vector $(\partial_y f)_{+}(x,0)$ is
a multiple of $e_3$ which proves (\ref{eqboundary4}). 
Lemma \ref{lemmaw12reflection} implies that the extension 
is of class $W^{2,2}(Q,\R^3)$. Noting that the $g_{ij}$ are even
functions, we conclude that $f$ is a $W^{2,2}$ conformal immersion
on $Q$ as defined in \cite{KL12}.\kasten

The next lemma states the linearized version of the constraints
(\ref{eqboundary1}), (\ref{eqboundary2}) on the boundary. Here 
$\eta$ denotes the inner (upward) normal of $Q_{+}$ along $I$ 
with respect to $g$, 
\begin{equation}
\label{eqconormal}
\eta = \frac{1}{\sqrt{g_{11} \det g}} \big(g_{11} e_2 - g_{12}e_1\big).
\end{equation} 

\begin{lemma} \label{lemmaadmissible} Let $f:(-\delta,\delta) \to 
W^{2,2}_{{\rm imm}}(Q_{+},\R^3)$,
$t \mapsto f(t)$, be a curve with derivative $\frac{d}{dt}f|_{t=0} = \phi$. If 
$f(t) \in {\cal M}$ for all $t \in (-\delta,\delta)$, i.e. (\ref{eqboundary1}) 
and (\ref{eqboundary2}) hold, then
\begin{eqnarray}
\label{eqadmissible1}
\langle \phi, e_3 \rangle & = & 0 \quad \mbox{ along }I,\\
\label{eqadmissible2}
\langle \partial_\eta \phi,\nu \rangle & = & 0 \quad \mbox{ along }I.
\end{eqnarray}
We say that $\phi \in W^{2,2} \cap W^{1,\infty}(Q_{+},\R^3)$ is
admissible at $f$ if  it satisfies  (\ref{eqadmissible1}), (\ref{eqadmissible2}).
\end{lemma}

{\em Proof. }(\ref{eqadmissible1}) is obvious. At $t=0$ the 
we have the following derivatives: 
\begin{eqnarray*}
g_{\alpha \beta}:(-\delta,\delta) \to W^{1,2} \cap L^\infty(Q_{+}), &
\displaystyle{\frac{d}{dt}}g_{\alpha \beta}|_{t=0} & 
= \langle \partial_\alpha \phi,\partial_\beta f \rangle 
+ \langle \partial_\alpha f,\partial_\beta \phi \rangle,\\
\nu:(-\delta,\delta) \to W^{1,2} \cap L^\infty(Q_{+},\R^3), & 
\displaystyle{\frac{d}{dt}}\nu|_{t=0} & 
= - g^{\alpha \beta} \langle \nu,\partial_\alpha \phi \rangle \partial_\beta f,\\
\eta:(-\delta,\delta) \to L^\infty(I,\R^2), & 
\displaystyle{\frac{d}{dt}}\eta|_{t=0} & = 
- \big(\langle \partial_\eta \phi,\partial_\tau f \rangle 
+ \langle \partial_\tau \phi,\partial_\eta f \rangle\big)\,\tau\\
&&\quad  - \langle \partial_\eta \phi,\partial_\eta f \rangle\,\eta.
\end{eqnarray*}
Here $\tau = e_1/\sqrt{g_{11}}$ is the unit tangent of $I$ with
respect to $g$. As $g$ is the metric induced by $f$, the vector 
$df \cdot \eta$ is by definition orthogonal to $df \cdot e_1$ 
and also to $\nu$. But these two vectors span $\R^2$ by (\ref{eqboundary1}) 
and (\ref{eqboundary2}), thus $df \cdot \eta = \pm e_3$. We now compute
at $t = 0$ along $I$ 
\begin{eqnarray*}
0 & = & \frac{d}{dt} \langle \nu,e_3 \rangle|_{t=0} 
\quad \quad \quad \mbox{ (using $\nu \perp e_3$ by (\ref{eqboundary2}))}\\
  & = & \pm \langle \frac{d}{dt}\nu|_{t=0},df \cdot \eta \rangle \quad
        \,\,\mbox{ (using $df \cdot \eta = \pm e_3$, see above)}\\
  & = & \mp \langle \nu,\frac{d}{dt}(df \cdot \eta)|_{t=0} \rangle  \quad
        \mbox{ (using $\nu \perp df \cdot \eta$)}\\
  & = & \mp \langle \nu, d\phi \cdot \eta \rangle \quad \quad
        \quad \quad \mbox{ (using $\nu \perp df \cdot \frac{d}{dt}\eta|_{t=0}$)}.
\end{eqnarray*}
In this calculation we actually used the traces of the functions 
$g_{ij}$, $\nu$, $\partial_i f$ and $\partial_i \phi$. The trace 
operator is continuous from $W^{1,2} \cap L^\infty(Q_{+})$ to 
$L^\infty(I)$, and therefore interchanges with the time derivative 
at $t = 0$. This justifies our computation.
\kasten

Next recall the first variation formula for the Willmore energy,
see \cite{Riv08}. The Willmore functional ${\cal W}(f)$ is 
Fr\'{e}chet differentiable on $W^{2,2}_{{\rm imm}}(Q_{+},\R^3)$, 
with derivative
\begin{eqnarray}
\nonumber
D{\cal W}(f)\phi & = & \frac{1}{2} \int_{Q_{+}} 
\langle \vec{H},\Delta_g \phi \rangle\,d\mu_g\\
\label{eqfirstvariation}
&& - \int_{Q_{+}} g^{ij}g^{kl} \langle \vec{H},A_{ik} \rangle 
\langle \partial_j f,\partial_k \phi \rangle\,d\mu_g\\
\nonumber
&& + \frac{1}{4} \int_{Q_{+}} |\vec{H}|^2 g^{ij} \langle \partial_i f,\partial_j \phi \rangle\,d\mu_g.
\end{eqnarray}

\begin{theorem} \label{thm1} Let $f \in W^{2,2}_{{\rm conf}}(Q_{+},\R^3)$ satisfy 
(\ref{eqboundary1}), (\ref{eqboundary2}) along $I$. Assume that $f$ is Willmore 
critical under these constraints, i.e. whenever 
$\phi \in W^{2,2}\cap W^{1,\infty}(Q_{+},\R^3)$ is admissible at $f$
with compact support in $Q_{+} \cup I$, then 
\begin{equation}
\label{eqlinecritical}
D{\cal W}(f)\phi = 0.
\end{equation}
Then extending $f$ by reflection at $\R^2$
yields a smooth Willmore immersion $f:Q \to \R^3$.
\end{theorem}

{\em Proof. }We establish the weak version of the Willmore equation on all of $Q$. 
Let $R = {\rm diag}(1,1,-1)$. We say that $\phi:Q \to \R^3$ is even (resp. odd) 
iff $\phi(x,y) = \pm R \phi(x,-y)$. By definition, the immersion $f$ and 
its derivatives $\partial_1 f,\partial_1^2 f,\partial_2^2 f$ are even, 
while $\partial_2 f,\partial^2_{12} f$ are odd. The derivative $\partial_1$ 
preserves the parity, while $\partial_2$ changes it. One checks that the 
integrand in the first variation formula is odd in the case when $\phi$ is 
odd, so that the integral vanishes trivially in this case. For $\phi$ even
the integrand is even and 
$$
D{\cal W}(f)\phi = 2\, D{\cal W}(f|Q_{+})\phi|_{Q_{+}}.
$$
As $\phi$ is even, we have $\phi_3(x,0) = 0$ and $\partial_y \phi_i(x,0) = 0$ 
for $i = 1,2$. Since $f|_{Q_{+}}$ is conformally parametrized, this
implies $\partial_\eta \phi_i (x,0) = 0$. Thus the variation $\phi$ 
is admissible, and by assumption 
$$
D{\cal W}(f|Q_{+})\phi|_{Q_{+}} = 0.
$$
It follows that the $W^{2,2}$ conformal immersion $f:Q \to \R^3$ 
is a critical point of the Willmore functional with respect to 
all compactly supported variations. The regularity theorem of 
Rivi\`{e}re \cite{Riv08} implies that $f$ is a smooth 
Willmore immersion. \kasten

\section{Nondegeneracy of the boundary condition}
\label{section2}
Here we address the question whether any admissible $\phi$ is the tangent
vector of some curve $f(t)$ of $W^{2,2}$ immersions satisfying the 
boundary constraints. Let $f \in W^{2,2}_{{\rm imm}}(Q_{+},\R^3)$ 
be a given conformal immersion
satisfying 
\begin{equation}
\label{eqconstraint}
{\cal B}(f) = \big(\langle f,e_3 \rangle,\langle \nu_f ,e_3 \rangle\big) = 0 
\quad \mbox{ along }I.
\end{equation}
Formally linearizing ${\cal B}$ at $f$ yields the operator 
\begin{equation}
\label{eqlinearized}
L_f \phi = \big(\langle \phi,e_3 \rangle,\langle D\phi \cdot \eta, \nu \rangle\big).
\end{equation}
For $f$ conformal, i.e. $g_{ij} = e^{2u} \delta_{ij}$,  
the operator $L_f$ becomes
\begin{equation}
\label{eqlinearizedconf}
L_f \phi = \big(\langle \phi,e_3 \rangle, e^{-u} \langle \partial_y \phi, \nu \rangle\big).
\end{equation}
Putting $K = [-\frac{\pi}{2},\frac{\pi}{2}]$  we consider the Banach spaces 
\begin{eqnarray*}
X & = & \{(a,b) \in W^{\frac{3}{2},2}\cap W^{1,\infty}(I) 
\oplus W^{\frac{1}{2},2}\cap L^\infty(I): {\rm spt\,}(a,b) \subset K\},\\
Y & = &  \{\phi \in W^{2,2} \cap W^{1,\infty}(Q_{+},\R^2 \oplus \R):
{\rm spt\,}L_f\phi \subset K\}.
\end{eqnarray*}
We have the linear map, using the extension $BH$ from the appendix,
\begin{equation}
\label{eqextension}
\Phi_f:X \to Y,\,\,\Phi_f(a,b) = BH(0,e^u b \nu_f) \oplus BH(a,0)e_3.
\end{equation}
Note that $\nu_f:I \to \R^2$ by (\ref{eqconstraint}). Since 
$u \in W^{1,2} \cap L^\infty(Q_{+})$ and $\nu_f \in W^{1,2} \cap L^\infty(Q_{+},\R^3)$ 
by assumption, we can estimate
\begin{equation}
\label{eqextensionbounds}
\|\Phi_f(a,b)\|_Y \leq C\,\|(a,b)\|_X. 
\end{equation}
The constant depends on the bound for $f$ in $W^{2,2} \cap W^{1,\infty}(Q_{+},\R^3)$, 
and on the lower bound for $u$, i.e. for the Jacobian of $f$. Now $L_f$
is a continuous operator
\begin{eqnarray*}
L_f:Y \to X,\,
L_f(\phi) = \big(\langle \phi,e_3 \rangle, e^{-u} \langle \partial_y \phi,\nu \rangle\big).
\end{eqnarray*}
Clearly $\phi \in {\rm ker\,}L_f$ if and only if 
$\langle \phi, e_3 \rangle = 0$, $\langle \partial_y \phi,\nu \rangle = 0$,
along $I$. By construction 
$$
L_f \circ \Phi_f = {\rm id}, \quad \mbox{ in particular } 
{\rm ker\,}L_f \cap {\rm im\,}\Phi_f = \{0\}.
$$ 
The subspace ${\rm im\,}\Phi_f \subset Y$ is closed:
let $\phi_k = \Phi_f(a_k,b_k) \to \phi$ in $Y$. 
Then $(a_k,b_k) = L_f \phi_k \to L_f \phi =: (a,b)$ in $X$,
which yields $\phi = \lim_{k \to \infty} \Phi_f(a_k,b_k) = \Phi_f(a,b)$. 
It follows that restricting $L_f$ gives the isomorphism 
$$
L_f|_{{\rm im\,}\Phi_f}:{\rm im\,}\Phi_f \to X, 
$$
with inverse bounded by (\ref{eqextensionbounds}). 
Moreover we have the direct sum decomposition 
$$
Y = {\rm im\,}\Phi_f \oplus {\rm ker\,}L_f, \quad
\phi = \Phi_f (L_f \phi) \oplus \big(\phi - \Phi_f (L_f \phi)\big).
$$
By the implicit function theorem there exist open 
neighborhoods $U \subset {\rm ker\,}L_f$ and $V \subset {\rm im\,}\Phi_f$
of the origin, and a $C^1$ mapping $G:U \to V$ with  $G(0) = 0$,
such that for $\phi \oplus \psi \in U \oplus V$ the following 
equivalence holds:
$$
f + \phi + \psi\, \mbox{ satisfies (\ref{eqconstraint})} \quad 
\Leftrightarrow \quad \psi = G[\phi].
$$
Here we use that $W^{2,2}_{{\rm imm}}(Q_{+},\R^3)$ is an open subset 
of $W^{2,2} \cap W^{1,\infty}(Q_{+},\R^3)$, whence for $U,V$ 
sufficiently small the map $f+\phi + \psi$ is weakly immersed. 
For any $\phi \in {\rm ker\,}L_f$ we thus 
obtain the admissible curve of weak immersions
$$
f(t) = f + t \phi + G[t\phi],\quad
t \in (-\delta,\delta).
$$
We have $DG[0] = 0$ from the general construction. To see this explicitely, 
we note that $DG[0]\phi \in {\rm im\,}\Phi_f$ by definition, while 
differentiating at $t = 0$ shows $DG[0]\phi \in {\rm ker\,}L_f$:
$$
0 = \frac{d}{dt}{\cal B}(f(t))|_{t=0} = L_f (\phi + DG[0]\phi) = L_f DG[0]\phi.
$$
 
\section{Reflection in a line}

The condition of prescribing a line for the boundary
with free tangent plane is similar and in fact 
simpler than the previous one. We indicate the
main points.

\label{section3}
\begin{lemma} \label{lemmaw22line} Let $f:Q_{+} \to \R^3$ be a $W^{2,2}$
conformal immersion, that is
\begin{eqnarray}
\label{eqlineimmersed}
|\partial_x f \times \partial_y f| & \geq & \mu > 0 \quad \mbox{ in }Q_{+},\\
\label{eqlineconf1} 
\langle \partial_x f,\partial_y f \rangle & = & 0 \quad \mbox{ in }Q_{+},\\
\label{eqlineconf2}
|\partial_x f|^2 - |\partial_y f|^2 & = & 0 \quad \mbox{ in }Q_{+}.
\end{eqnarray}
Assume that $f|_I$ maps into the line $L = \{0\} \times \R \subset \R^3$, i.e.
\begin{equation}
\label{eqlineboundary}
f_i  = 0 \quad \mbox{ along $I$ for }i = 1,2. 
\end{equation}
Let $f$ be extended to all of $Q$ by reflection at $L$, that is
\begin{eqnarray}
\label{eqlineextend1}
f_i(x,y) & = &\,\, -f_i(x,-y) \quad \mbox{ for } i = 1,2,\\
\label{eqlineextend2}
f_3(x,y) & = & f_3(x,-y). 
\end{eqnarray}
Then $f:Q \to \R^3$ is a $W^{2,2}$ conformal immersion.
\end{lemma}

{\em Proof. }The function $f_3$ is even, while $f_i$, $i =1,2$, are 
odd with $f_i(x,0) = 0$ by assumption (\ref{eqboundary1}). Thus 
$f \in W^{1,2}(Q,\R^3)$ by Lemma \ref{lemmaw12reflection}. To see 
that $f \in W^{2,2}(Q,\R^3)$ we need to show the vanishing
of the odd first derivatives, i.e.
$$
\partial_x f_i = 0 \mbox{ on }I \mbox{ for } i=1,2  \quad \mbox{ and } \quad
\partial_y f_3 = 0 \mbox{ on } I.
$$
The first follows by differentiating the equation $f_i(x,0) = 0$ (copy the 
argument for $f_3$ from Lemma 2.2). To show the second statement at a 
point $(x,0)$, we may assume $\partial_y f(x,0) \neq 0$. Then 
$\partial_x f(x,0)$ is also nonzero. But $\partial_x f(x,0)$
is a multiple of $e_3$. By conformality, $\partial_y f(x,0)$ 
then lies in $\R^2$, i.e. $\partial_y f_3(x,0) = 0$. \kasten

\begin{theorem} \label{thm2} Let $f:Q_{+} \to \R^3$ be a $W^{2,2}$ conformal 
immersion with $f_i = 0$ for $i = 1,2$ along $I$.
Assume that $f$ is critical under these constraints, i.e. 
\begin{equation}
\label{eqcritical}
D{\cal W}(f)\phi = 0 \quad \mbox{ for all admissible }
\phi \in  W^{2,2} \cap W^{1,\infty}(Q_{+},\R^3),
\end{equation}
that is $\phi$ has compact support in $Q_{+} \cup I$ and 
$\phi_i|_I = 0$ for $i = 1,2$. Then extending $f$ by reflection 
at $L = \{0\} \times \R$ yields a smooth Willmore immersion 
$f:Q \to \R^3$.
\end{theorem}

{\em Proof. }We establish the weak version of the 
Willmore equation on all of $Q$.
Let $S = {\rm diag}(-1,-1,1)$. We say that $\phi:Q \to \R^3$ is even (resp. odd)
iff $\phi(x,y) = \pm S \phi(x,-y)$. By definition, the immersion $f$ and
its derivatives $\partial_1 f, \partial^2_{11} f, \partial^2_{22} f$ are even,
while $\partial_2 f,\partial_{12} f$ are odd. The derivative $\partial_1$
preserves the parity, while $\partial_2$ changes it. One checks {\bf (!)} that 
the integrand in the first variation formula is odd in the case when $\phi$ is
odd, so that the integral vanishes trivially in this case. For $\phi$ even
the integrand is even and one obtains
$$
D{\cal W}(f)\phi = 2\, D{\cal W}(f|Q_{+})\phi|_{Q_{+}}.
$$
Now $\phi$ even implies $\phi_i(x,0) = 0$ for $i = 1,2$, and hence 
the variation $\phi$ is admissible on $Q_{+}$, so that by assumption
$$
D{\cal W}(f|Q_{+})\phi|_{Q_{+}} = 0.
$$
It follows that the $W^{2,2}$ conformal immersion $\phi:Q \to \R^3$
is a critical point of the Willmore functional. The regularity
theorem of Rivi\`{e}re \cite{Riv08} implies that $f$ is a smooth
Willmore immersion. \kasten

In contrast to section 3, here both boundary constraints are
linear, hence the affine variation suffices 
to generate any given admissible field $\phi$.

\section{Critical points of related curvature energies}
\label{thomsen}

In this section we discuss the same subject for  
the two related functionals
\begin{eqnarray} 
\label{eqL2curvature}
{\mathcal E}(f) & = & \frac{1}{2} \int_{\Sigma} |h|^2\,d\mu_g,\\
\label{eqthomsen}
{\mathcal T}(f) & = & \frac{1}{2} \int_{\Sigma} |h^\circ|^2\,d\mu_g.
\end{eqnarray}
As is well-known, the integrand of the Thomsen functional 
${\mathcal T}(f)$ is pointwise conformally invariant \cite{Tho23}. 
Clearly $|h|^2 = |h^\circ|^2 + \frac{1}{2}H^2$, hence 
\begin{equation}
\label{eqfunctionals}
{\mathcal E}(f) = {\mathcal T}(f) + {\mathcal W}(f).
\end{equation} 
Assuming enough regularity we may rewrite the functionals using 
Gau{\ss}-Bonnet, i.e. 
\begin{eqnarray}
\label{eqL2gaussbonnet}
{\mathcal E}(f) & = & 2\, {\mathcal W}(f) 
+ \int_{\partial \Sigma} \varkappa_g\,ds - 2\pi \chi(\Sigma),\\
\label{eqthomsengaussbonnet}
{\mathcal T}(f) & = & {\mathcal W}(f)
+ \int_{\partial \Sigma} \varkappa_g\,ds - 2\pi \chi(\Sigma).
\end{eqnarray}
Therefore both functionals have $W(f)$ as Euler-Lagrange 
operator,  with factor $\frac{1}{2}$ for ${\mathcal T}(f)$.
Let $N^S$ be the interior unit normal along $S = \partial \Omega$,
and let $\eta$ be the interior unit normal along $\partial \Sigma$
with respect to the induced metric $g$. As before we 
denote by ${\widetilde{\mathcal M}}(S)$ the class of smooth immersions
$f:\Sigma \to \R^3$ satisfying the constraints
\begin{equation}
\label{eqconstraints}
f(\partial \Sigma) \subset S \quad \mbox{ and } \quad
\frac{\partial f}{\partial \eta} = N^S \circ f.
\end{equation}
Differentiating in the direction of the unit tangent 
$\tau = \frac{\partial}{\partial s}$ along $\partial \Sigma$ 
we obtain, using also $\langle \nu, N^S \circ f \rangle = 0$,
$$
0 = \partial_s \langle \partial_s f, N^S \circ f \rangle 
= \big \langle Df \cdot \nabla_\tau \tau, N^S \circ f \big\rangle 
+ \big\langle \partial_s f, (DN^S) \cdot \partial_s f \big\rangle.
$$
Now by definition $\nabla_\tau \tau = \varkappa_g \eta$. 
Using again (\ref{eqconstraints}) we see that 
\begin{equation}
\label{eqgeodesiccurvature}
\varkappa_g = h^S \big(\partial_s f, \partial_s f\big). 
\end{equation}
We conclude that for $f \in {\widetilde{\cal M}}(S)$ we have
\begin{eqnarray}
\label{eqL2gaussbonnet}
{\mathcal E}(f) & = & 2\, {\mathcal W}(f) + \int_{\partial \Sigma}
h^S \big(\partial_s f, \partial_s f\big)\,ds - 2\pi \chi(\Sigma),\\
\label{eqthomsengaussbonnet}
{\mathcal T}(f) & = & {\mathcal W}(f) + \int_{\partial \Sigma}
h^S \big(\partial_s f, \partial_s f\big)\,ds - 2\pi \chi(\Sigma).
\end{eqnarray}
On ${\widetilde{\mathcal M}}(\R^2)$ the three functionals $\mathcal{E}$, $\mathcal{T}$ and $\mathcal{W}$ coincide 
up to a topological constant, in particular they have the 
same critcial points. Next we compute the resulting free boundary conditions.

\begin{lemma} \label{lemmafirstvariation} For a smooth variation 
$f(\cdot,t):\Sigma \to \R^3$, $t \in (-\delta,\delta)$ with velocity 
$\phi = \varphi \nu + Df \cdot \xi$ we have 
\begin{eqnarray}
\label{eqL2variation}
\frac{d}{dt} {\mathcal E}(f(\cdot,t))|_{t=0} & = & 
\int_{\Sigma} {W}(f) \varphi \, d\mu_g +\int_{\partial \Sigma} 
 \tau (\eta)\,ds_g,
\quad \mbox{ where }\\
\tau(\eta) & = &  \frac{\partial H}{\partial \eta} \varphi  
-   h({\rm grad}_g \varphi,\eta)
- \frac{1}{2} |h|^2 g(\xi,\eta).
\end{eqnarray}
\end{lemma}

{\em Proof. }Assume first that $\phi$ is normal. The following 
identities in general codimension are computed in \cite{KS02}:
\begin{eqnarray*}
\partial_t g_{\alpha \beta} & = & - 2 h_{\alpha \beta}\varphi ,\\
\partial_t (d\mu_g) & = & - H\varphi \, d\mu_g,\\
\partial_t h_{\alpha \beta} & = & 
\nabla^2_{\alpha \beta}\varphi 
- g^{\lambda \mu} h_{\alpha \lambda}  h_{\beta \mu} \varphi .
\end{eqnarray*} 
We compute further, using normal coordinates at $t= 0$, $p \in \Sigma$,
\begin{eqnarray*}
\partial_t \Big(\frac{1}{2} |h|^2 d\mu_g\Big) & = & 
\frac12 \partial_t \big(g^{\alpha \lambda}g^{\beta \mu}
h_{\alpha \beta}h_{\lambda \mu}  \,d\mu_g\big)\\
& = & \big( \nabla^2_{\alpha \beta} \varphi 
- h_{\alpha \gamma} h_{\beta \gamma} \varphi \big)h_{\alpha \beta} \,d\mu_g
- \frac{1}{2} |h|^2\, H\varphi\,d\mu_g\\
&& +2\, h_{\alpha \lambda}  
h_{\alpha \beta}h_{\lambda \beta} \varphi \,d\mu_g.
\end{eqnarray*} 
Decomposing 
$h_{\alpha \beta} = h^\circ_{\alpha \beta} + \frac{1}{2} H \delta_{\alpha \beta}$
yields
\begin{equation}
\label{eqcubiccurvature}
h_{\alpha \lambda}   h_{\alpha \beta}h_{\lambda \beta}
= \frac32 |h^\circ|^2 H+\frac14 H^2H.
\end{equation}
We finally arrive at
\begin{equation}
\partial_t \Big(\frac{1}{2} |h|^2 d\mu_g\Big) = 
\big(h \nabla^2 \varphi  + |h^\circ|^2H\varphi \big)\,d\mu_g. 
\end{equation}
Now we compute, still in normal coordinates, using that 
$\nabla_i h_{ik} =  \nabla_k H$ by Codazzi,
\begin{eqnarray*}
h\nabla^2 \varphi  -  \varphi \Delta H 
& = &  \frac{1}{\sqrt{\det g}} 
\partial_\alpha \big(\sqrt{\det g}\, g^{\alpha \beta} \tau_\beta \big), 
\quad \mbox{ where }\\
\tau_\beta & = & g^{\lambda \mu}  \nabla_\lambda \varphi h_{\beta \mu}  
-  \varphi \nabla_\beta H .
\end{eqnarray*}
Next, consider a tangential variation of the form $f \circ \varphi_t$, where
$\varphi_t$ is the flow of a vectorfield $\xi$. Then we have
$$
{\mathcal E}(f \circ \varphi_t,\Omega) = {\mathcal E}(f,\varphi_t(\Omega))\\
= \frac{1}{2} \int_\Omega |h|^2(\varphi_t(x)) J_g\varphi_t(x)\,d\mu_g(x). 
$$
Differentiating at $t = 0$ we obtain
$$
\frac{d}{dt}{\mathcal E}(f \circ \varphi_t,\Omega)|_{t=0} = 
\frac{1}{2} \int_{\Omega} \big(d|h|^2(\xi) + |h|^2 {\rm div_g}\xi\big)\,d\mu_g
= \frac{1}{2} \int_{\Omega} {\rm div_g} (|h|^2 \xi)\,d\mu_g
$$
Now $g(|h|^2 \xi,X) = |h|^2 \langle \phi, df \cdot X \rangle$.
The claim of the lemma follows by combinig the two computations. 
\kasten

Now assume that $f$ is a critical point in ${\mathcal M}(S)$. As 
computed in \cite{AK14}, a variation $\phi = \varphi \nu + Df \cdot \xi$
is admissible if and only if 
\begin{equation}
\label{eqadmissibleclassic}
g(\xi,\eta) = 0 \quad \mbox{ and }\quad 
\frac{\partial \varphi}{\partial \eta} + \varphi h^S(\nu,\nu) = 0
\quad \mbox{ on } \partial \Sigma.
\end{equation}
It is easy to see that for given functions $\varphi,\mu$ on 
$\partial \Sigma$, there exists an admissible variation $\phi$
such that $\phi = \varphi \nu + \mu \partial_s f$ along 
$\partial \Sigma$. Now if $f$ is critical in ${\widetilde{\mathcal M}}(S)$,
then $W(f) = 0$ and further for $\phi$ admissible
$$
0 = \int_{\partial \Sigma} 
\Big(\varphi \frac{\partial H}{\partial \eta}
- h({\rm grad}_g \varphi,\eta) 
-\frac{1}{2}|h|^2 g(\xi,\eta)\Big)\,ds_g.
$$
We have $g(\xi,\eta) = 0$ from (\ref{eqadmissibleclassic}). Furthermore
\begin{eqnarray*}
h({\rm grad}_g \varphi,\eta) & = & 
(\partial_\tau \varphi) h(\tau,\eta) + (\partial_\eta \varphi) h(\eta,\eta)\\
& = & \partial_\tau (\varphi h(\tau,\eta)) - \varphi \nabla_\tau h(\tau,\eta)
- \varphi h(\nabla_\tau \tau,\eta) - \varphi h(\tau,\nabla_\tau \eta)\\
& = & \partial_\tau (\varphi h(\tau,\eta)) - \varphi \nabla_\tau h(\tau,\eta)
- \varphi \varkappa_g h(\eta,\eta) + \varphi \varkappa_g  h(\tau,\tau).
\end{eqnarray*}
We arrive at the formula
$$
0 = \int_{\partial \Sigma}
\Big(\varphi \frac{\partial H}{\partial \eta} - \frac{\partial \varphi}{\partial \eta} h(\eta,\eta)
+ \varphi\big[\nabla_\tau h(\tau,\eta) 
+ \varkappa_g \big(h(\eta,\eta)-h(\tau,\tau)\big)\big]\Big)\,ds_g.
$$
Now differentiating the constraint $\langle \nu,N^S \circ f \rangle =0$ in 
direction of $\tau$, we get
\begin{eqnarray*}
0 & = & \partial_\tau \langle \nu,N^S \circ f \rangle\\ 
& = &  \langle Df \cdot W \tau, Df \cdot \eta \rangle 
+ \langle \nu, (DN^S) \circ f Df \cdot \tau \rangle\\
& = & h(\tau,\eta) + h^S(\nu,Df \cdot \tau).
\end{eqnarray*}
Differentiating once more yields
$$
\nabla_\tau h(\tau,\eta) + \varkappa_g \big(h(\eta,\eta)-h(\tau,\tau)\big)
+ \partial_s \big[h^S(\nu, \partial_s f)\big] = 0.
$$
Inserting this into the boundary condition, we finally get
\begin{eqnarray*}
0 & = & \int_{\partial \Sigma} 
\Big(\varphi \frac{\partial H}{\partial \eta} - \frac{\partial \varphi}{\partial \eta} h(\eta,\eta)
- \varphi\,\partial_s \big[h^S(\nu, \partial_s f)\big]\Big)\,ds_g\\ 
& = & \int_{\partial \Sigma} \varphi\,
\Big(\frac{\partial H}{\partial \eta} + h(\eta,\eta)h^S(\nu,\nu)
- \partial_s \big[h^S(\nu, \partial_s f)\big]\Big)\,ds_g.
\end{eqnarray*}
Thus we obtained the boundary condition
\begin{equation}
\label{eqL2boundarycondition}
\frac{\partial H}{\partial \eta} + h^S(\nu,\nu)h(\eta,\eta)
- \partial_s \big[h^S(\nu, \partial_s f)\big] = 0.
\end{equation}
Clearly this reduces to $\frac{\partial H}{\partial \eta} = 0$ 
when $S$ is a plane. For a sphere $S = \partial B_R(0)$ we get
$$
\frac{\partial H}{\partial \eta} + \frac{1}{R} h(\eta,\eta) = 0.
$$
To get the boundary condition for the Thomsen functional $\mathcal{T}$ one
just combines the results for ${\mathcal W}(f)$ and 
${\mathcal E}(f)$; we just state the equations:
\begin{eqnarray*}
0 & = & \frac{\partial H}{\partial \eta} + h^S(\nu,\nu) H = 0,\\
0 & = & \frac{\partial H}{\partial \eta} 
+ h^S(\nu,\nu) (h(\eta,\eta-h(\tau,\tau))
- \partial_s \big[h^S(\nu, \partial_s f)\big]. 
\end{eqnarray*}
The above computations all assumed sufficient regularity. For the application of our main theorem we need \eqref{eqL2gaussbonnet} and \eqref{eqthomsengaussbonnet} for conformal immersions which are only in $\mathcal{M}$. This follows directly from a reflection argument using Lemma \ref{lemmaw22conformal} and the Gauss-Bonnet theorem for closed surfaces from \cite{KL12}.  \\
\\

\section{Appendix: Extension lemma}
\label{appendixb}

Here we extend compactly supported functions on $\R$
to the upper halfplane with certain bounds. Such an extension 
has been also mentioned in \cite{Sim93}. In addition 
to the previous notation we put $K = [-\frac{\pi}{2},\frac{\pi}{2}]$
and $\H = \R \times [0,\infty)$.

\begin{lemma}[Extension] 
Let $\varphi \in W^{\frac{3}{2},2} \cap W^{1,\infty}(\R)$,
$\psi \in W^{\frac{1}{2},2} \cap L^\infty(\R)$ have support 
in $K$. There exists
$u \in W^{2,2} \cap W^{1,\infty}(\H)$ with
compact support in $I \times [0,1)$, such that
\begin{equation}
\label{eqboundarydata}
u(\cdot,0) = \varphi \quad \mbox{ and } \quad 
\frac{\partial u}{\partial y}(\cdot,0) = \psi,
\end{equation}
and such that for a universal constant $C< \infty$ one 
has the estimates 
\begin{eqnarray}
\label{eqc1estimate}
\|u\|_{C^1} & \leq &
C \,\big(\|\varphi\|_{W^{1,\infty}} + \|\psi\|_{L^\infty}\big),\\
\label{eql2estimate}
\|u\|_{W^{2,2}} & \leq &
C \, \big(\|\varphi\|_{W^{\frac{3}{2},2}} + \|\psi\|_{W^{\frac{1}{2},2}}\big),
\end{eqnarray}
\end{lemma} 

{\em Proof. }We start by considering the harmonic extension 
\begin{equation}
\label{eqextharmonic}
H\varphi:\H \to \R,\,
H\varphi(x,y) = \sum_{k=0}^\infty \varphi_k(x) e^{-ky}.
\end{equation}
Here $\varphi = \sum_{k=0}^\infty \varphi_k$ is the 
Fourier decomposition on $I$. We compute 
\begin{eqnarray*}
H\varphi(x,y) & = &
\frac{1}{\pi}\sum_{k=0}^\infty
\Big(\int_I \varphi(x') \cos (kx')\,dx'\Big) \cos kx\, e^{-ky}\\
&& + \frac{1}{\pi}\sum_{k=0}^\infty
\Big(\int_I \varphi(x') \sin (kx')\,dx'\Big) \sin kx\, e^{-ky}\\
& = & \frac{1}{\pi} \int_I \varphi(x') \sum_{k=0}^\infty \cos k(x-x')\, e^{-ky}\,dx'\\
& = & \frac{1}{\pi}\, {\rm Re\,}\int_I \varphi(x') \sum_{k=0}^\infty e^{k(-y+i(x-x'))}\,dx'\\
& = & \frac{1}{\pi}\, {\rm Re\,}\int_I \varphi(x')\frac{1}{1- e^{-y+i(x-x')}}\,dx'\\
& = & \frac{1}{\pi}\, \int_I \varphi(x')
\frac{1-e^{-y}\cos (x-x')}{1-2e^{-y} \cos (x-x') + e^{-2y}}\,dx'\\
& = & \frac{1}{2\pi}\, \int_I \varphi(x') 
\frac{e^y - \cos (x-x')}{\cosh y - \cos (x-x')}\,dx'.
\end{eqnarray*}
Thus we have 
$$
H\varphi(x,y) = \int_I G(x-x',y)\varphi(x')\,dx' \quad 
\mbox{ where }G(x,y) = \frac{1}{2\pi}\,\frac{e^y - \cos x}{\cosh y - \cos x}.
$$
For $\varphi(x) \equiv 1$ we have $H \varphi \equiv 1$, which yields
noting $G(-x,y) = G(x,y)$,
$$
\int_I G(x,y)\,dx = 1 \quad \mbox{ for all }y > 0.
$$
In particular for any $y > 0$ we can estimate
\begin{equation}
\label{eqharmonic1}
\|H\varphi(\cdot,y)\|_{C^0(I)} \leq \|\varphi\|_{L^\infty(I)}.
\end{equation}
Using the Cauchy estimate on the disk $D_y((x,y))$ we get, again for $y> 0$, 
\begin{equation}
\label{eqharmonic2}
\|D(H \varphi)(\cdot,y)\|_{C^0(I)} 
\leq \frac{C}{y} \,\|\varphi\|_{L^\infty(I)}.
\end{equation}
Now the biharmonic extension $u = BH(\varphi,\psi)$ is given by 
$$
u(x,y) = H\varphi(x,y) - y\, \partial_y(H\varphi)(x,y) + y H\psi(x,y).
$$
The extension has initial values 
$$
u(x,0) = H\varphi(x,0) = \varphi(x) \quad \mbox{ and } \quad
\partial_y u(x,0) = H\psi(x,0) = \psi(x).
$$
To see that $u(x,y)$ is biharmonic, we note for $h(x,y)$ harmonic that 
$$
\Delta (y h(x,y)) = 2\, \partial_y h(x,y)  \quad \mbox{ thus }\quad
\Delta^2 (y h(x,y)) = 2 \partial_y \Delta h(x,y) = 0. 
$$
Now we collect the relevant estimates. By (\ref{eqharmonic1}) and
(\ref{eqharmonic2}) we get for $0 < y \leq 1$ 
\begin{equation}
\label{eqbiharmonic1}
\|u(\cdot,y)\|_{C^0(I)} \leq C \|\varphi\|_{L^\infty(I)} + y \|\psi\|_{L^\infty(I)}.
\end{equation}
Writing $\partial_x u = H\varphi' -y \,\partial_y (H\varphi') + y \partial_x (H\psi)$
we get by the above estimates 
\begin{equation}
\label{eqbiharmonic2}
\|\partial_x u(\cdot,y)\|_{C^0(I)} \leq 
C \big(\|\varphi'\|_{L^\infty(I)} + \|\psi\|_{L^\infty(I)}\big).
\end{equation}
Now we have using that $H \varphi$ is harmonic
$$
\partial_y u = - y\, \partial_y^2 (H\varphi) + H\psi + y\, \partial_y(H\psi)
= y\, \partial_x (H\varphi') + H\psi + y\, \partial_y(H\psi),
$$
whence by the above estimates 
\begin{equation}
\label{eqbiharmonic3}
\|\partial_y u(\cdot,y)\|_{C^0(I)} \leq 
C\, \big(\|\varphi'\|_{L^\infty(I)} + \|\psi\|_{L^\infty(I)}\big).
\end{equation}
For the $L^2$ estimates, we have for any $s \in \R$ by 
orthogonality of the $\varphi_k$
\begin{eqnarray*}
\big\|\sum_{k=1}^\infty k^s \varphi_k(x) e^{-ky}\|^2_{L^2(I \times (0,\infty))}
& = & \sum_{k = 1}^\infty k^{2s} \|\varphi_k\|_{L^2(I)}^2\,\int_0^\infty e^{-2ky}\,dy\\
& = & \frac{1}{2} \sum_{k=1}^\infty k^{2s-1} \|\varphi_k\|_{L^2(I)}^2\\
& = & \frac{1}{2} \big[\varphi\big]^2_{W^{s - \frac{1}{2},2}(I)}.
\end{eqnarray*}
Now for $i,j \in \N_0$ we have 
$\partial_x^i \partial_y^j\, (H\varphi) 
= (-1)^j \sum_{k=0}^\infty \varphi_k^{(i)}(x)\,k^je^{-ky}$.
As $\|\varphi_k^{(i)}\|_{L^2(I)} = k^i\,\|\varphi_k\|_{L^2(I)}$, we 
get by putting $s = 0,1,2$ 
\begin{eqnarray*}
\|H(\varphi-\varphi_0)\|_{L^2(I \times (0,\infty))} 
& \leq & C [\varphi]_{W^{-\frac{1}{2},2}(I)},\\
\|D(H \varphi)\|_{L^2(I \times (0,\infty))} & \leq & C [\varphi]_{W^{\frac{1}{2},2}(I)},\\
\|D^2(H \varphi)\|_{L^2(I \times (0,\infty))} & \leq & C [\varphi]_{W^{\frac{3}{2},2}(I)}.
\end{eqnarray*}
In the following we put $\|\varphi\|^2_{W^{s,2}(I)} = |\varphi_0|^2 + [\varphi]^2_{W^{s,2}(I)}$.
For the biharmonic extension $u$ we deduce the following estimates: 
\begin{eqnarray*}
\|u\|_{L^2(I \times (0,1))} & \leq & 
C\, \big(\|\varphi\|_{W^{\frac{1}{2},2}(I)} + \|\psi\|_{W^{-\frac{1}{2},2}(I)}\big),\\
\|D u\|_{L^2(I \times (0,1))} & \leq &
C\, \big([\varphi]_{W^{\frac{3}{2},2}(I)} + [\psi]_{W^{\frac{1}{2},2}(I)}\big),\\
\|\Delta u\|_{L^2(I \times (0,1))} & = & 
2\, \|-\partial_y^2 (H\varphi) + \partial_y (H\psi)\|_{L^2(I \times (0,1))}\\ 
& = & 2\, \|H\varphi'' + \partial_y (H\psi)\|_{L^2(I \times (0,1))}\\
& \leq & C\, \big([\varphi'']_{W^{-\frac{1}{2},2}(I)} + [\psi]_{W^{\frac{1}{2},2}(I)}\big)\\
& = & C\, \big([\varphi]_{W^{\frac{3}{2},2}(I)} + [\psi]_{W^{\frac{1}{2},2}(I)}\big).
\end{eqnarray*}
To estimate the full second derivatives we integrate by parts:
\begin{eqnarray*}
\int_{I \times (0,1)} |D^2 u|^2 & = & \int_{I \times (0,1)} |\Delta u|^2 
+ \int_{I \times (0,1)} \sum_{i,j = 1}^2 \Big(
\partial_i \big(\partial_j u\, \partial^2_{ij} u \big)
-\partial_j \big(\partial_j u\, \partial_{ii}^2 u\big)\Big)\\
& = & \int_{I \times (0,1)} |\Delta u|^2  +\left[\int_{I \times \{y\}} 
\big(\partial_j u\, \partial^2_{2j} u 
- \partial_2 u\, \partial^2_{ii} u\big)\,dx\right]_{y=0}^{y = 1}\\
& = & \int_{I \times (0,1)} |\Delta u|^2  +\left[\int_{I \times \{y\}}
\big(\partial_x u\, \partial^2_{xy} u - \partial_y u\, \partial^2_{xx} u\big)\,dx\right]_{y=0}^{y = 1}.
\end{eqnarray*}
At $y = 0$ we get
$$
\int_{I \times \{0\}}
\big(\partial_x u\, \partial^2_{xy} u - \partial_y u\, \partial^2_{xx} u\big)\,dx
= \int_{I} (\varphi' \psi' - \varphi'' \psi)\,dx
= -2 \sum_{k =1}^\infty k^2 \langle \varphi_k, \psi_k \rangle_{L^2(I)}.
$$
The Cauchy-Schwarz iequality yields 
\begin{eqnarray*}
\Big|\int_{I \times \{0\}}
\big(\partial_x u\, \partial^2_{xy} u - \partial_y u\, \partial^2_{xx} u\big)\,dx\Big|
& \leq & 2 \Big(\sum_{k=1}^\infty k^3 \|\varphi_k\|_{L^2(I)}^2 \Big)^{\frac{1}{2}}
\Big(\sum_{k=1}^\infty k \|\psi_k\|_{L^2(I)}^2 \Big)^{\frac{1}{2}}\\
& = & 2 [\varphi]_{W^{\frac{3}{2},2}(I)} [\psi]_{W^{\frac{1}{2},2}(I)}. 
\end{eqnarray*}
On the other hand, standard interior estimates imply that 
$$
\Big|\int_{I \times \{1\}}
\big(\partial_x u\, \partial^2_{xy} u - \partial_y u\, \partial^2_{xx} u\big)\,dx\Big|
\leq C\,(\|\varphi\|_{L^2(I)} + \|\psi\|_{L^2(I)}).
$$
We have proved that
\begin{eqnarray}
\label{eqc1estimate}
\|u\|_{C^1(I \times (0,1))} & \leq &
C \,\big(\|\varphi\|_{W^{1,\infty}(I)} + \|\psi\|_{L^\infty(I)}\big),\\
\label{eql2estimate}
\|u\|_{L^2(I \times (0,1))} & \leq &
C\,\big(\|\varphi\|_{W^{\frac{1}{2},2}(I)} + \|\psi\|_{W^{-\frac{1}{2},2}(I)}\big),\\
\label{eq1stl2estimate}
\|D u\|_{L^2(I \times (0,1))} & \leq &
C\, \big([\varphi]_{W^{\frac{3}{2},2}(I)} + [\psi]_{W^{\frac{1}{2},2}(I)}\big),\\
\label{eq2ndl2estimate}
\|D^2 u\|_{L^2(I \times (0,1))} & \leq & 
C \, \big([\varphi]_{W^{\frac{3}{2},2}(I)} + [\psi]_{W^{\frac{1}{2},2}(I)}\big),
\end{eqnarray}
The estimate for $\|Du\|_{L^2}$ is not optimal. To conclude the proof 
of the lemma we choose a cutoff function $\eta \in C_c^\infty(I \times [0,1))$
such that
$$
\eta \equiv 1 \quad \mbox{ on } \quad  
\Big[-\frac{\pi}{2}, \frac{\pi}{2}\Big] \times \Big[0,\frac{1}{2}\Big].
$$
The function $\eta u$ has all properties required in the lemma. \kasten

\end{document}